\documentclass{amsart}

\usepackage[all]{xy}

\usepackage[a4paper, top=2cm, bottom=2cm, left=3cm, right=3cm, heightrounded, marginparwidth=3.5cm, marginparsep=.2cm, centering]{geometry}
\usepackage{hyperref}
\usepackage{cleveref}

\usepackage{bbold}
\usepackage{amssymb}
\usepackage{amsaddr}
\usepackage{todonotes}
\usepackage[T1]{fontenc}
\usepackage[backend=bibtex, style=numeric]{biblatex}

\theoremstyle{plain}
\newtheorem{theorem}[subsection]{Theorem}
\newtheorem*{theorem*}{Theorem}
\newtheorem{proposition}[subsection]{Proposition}
\newtheorem{lemma}[subsection]{Lemma}
\newtheorem{corollary}[subsection]{Corollary}
\newtheorem{conjecture}[subsection]{Conjecture}
\theoremstyle{definition}
\newtheorem{definition}[subsection]{Definition}

\newcommand{\Bun}{\mathrm{Bun}}
\newcommand{\calO}{\mathcal{O}}
\newcommand{\cont}{\mathrm{cont}}

\newcommand{\ex}{\mathrm{ex}}
\newcommand{\F}{\mathbb{F}}

\newcommand{\FF}{\mathrm{FF}}

\newcommand{\Fun}{\mathrm{Fun}}

\newcommand{\Hom}{\mathrm{Hom}}

\newcommand{\Isoc}{\mathrm{Isoc}}

\newcommand{\Perf}{\mathrm{Perf}}
\newcommand{\proet}{\mathrm{pro\acute{e}t}}
\newcommand{\Q}{\mathbb{Q}}
\newcommand{\Rep}{\mathrm{Rep}}

\newcommand{\sch}{\mathrm{sch}}
\newcommand{\Sets}{\mathrm{Sets}}
\newcommand{\Spa}{\mathrm{Spa}}
\newcommand{\Spd}{\mathrm{Spd}}
\newcommand{\Spec}{\mathrm{Spec}}

\newcommand{\Z}{\mathbb{Z}}
\newcommand{\R}{\mathbb{R}}

\newcommand{\cal}[1]{\mathcal{#1}}

\DeclareMathOperator{\GL}{GL}

\DeclareMathOperator{\pr}{pr}

\begin{document}

\title{$G$-bundles on the absolute Fargues--Fontaine curve}
\author{Johannes Ansch\"{u}tz}
\date{02.06.2016}

\begin{abstract}
  We prove that the category of ``vector bundles on the absolute Fargues--Fontaine curve'' (more precisely the category of sections over some discrete algebraically closed field of the $v$-stack $\Bun_\FF$ of vector bundles on the Fargues--Fontaine curve) is canonically equivalent to the category of isocrystals. We deduce a similar result for ``$G$-bundles on the absolute Fargues--Fontaine curve'' for some reductive group $G$ as well as for sections of $\Bun_\FF$ over classifying stacks for locally profinite groups.
\end{abstract}

\maketitle

\section{Introduction}
\label{sec:introduction-}

Let $p$ be a prime, and let $E$ be a non-archimedean local field with residue field $\F_q$ of characteristic $p$. Let $\Perf_{\F_q}$ be the category of perfectoid spaces over $\F_q$, and let $\Bun_{\FF}$ be the $v$-stack on $\Perf_{\F_q}$ of vector bundles on the Fargues--Fontaine curve for $E$, cf.\ \cite[Section II.2]{fargues2021geometrization}.
Thus, if $S$ is any perfectoid space over $\F_q$, then $\Bun(S)$ is the category of vector bundles on the (adic) relative Fargues--Fontaine curve $X_{E,S}$, cf.\ \cite[Definition II.1.15]{fargues2021geometrization}. 
Given now any small $v$-stack $Y$ over $\F_q$ (in the sense of \cite[Definition 12.4]{scholze_etale_cohomology_of_diamonds}) one can contemplate the category
\[
  \Bun_\FF(Y)
\]
of morphisms of $v$-stacks $Y\to \Bun_\FF$ over $\Perf_{\F_q}$. More concretely, if $Y$ is the quotient of a perfectoid space $X$ by some equivalence relation $R=X\times_Y X$, which is again represented by some perfectoid space, then $\Bun_\FF(Y)$ is equivalent to the category of descent data for $X\to Y$. As the pullback of vector bundles is exact, $\Bun_\FF(Y)$ aquires a natural structure of an exact category.

The category $\Bun_\FF(Y)$ for small $v$-sheaves associated with adic spaces over $\Z_p$ occur naturally when contemplating integral models of local Shimura varieties, cf.\ \cite{pappas2021p}. It seems therefore to be an interesting question to understand the category $\Bun_\FF(Y)$ for some specific examples of small $v$-sheaves $Y$, e.g., the small $v$-sheaf
\[
\Spd(R):=\Spec(R)^\diamond\colon \Perf_{\F_q}\to (\Sets),\ S\mapsto \Hom_{\F_q}(R,\mathcal{O}_S(S))
\]
associated with some perfect $\F_q$-algebra $R$, cf.\ \cite[Section 18.3]{scholze2020berkeley}.

The main result of this paper discusses the case that $Y=\Spd(k)$ is the small $v$-sheaf associated with a (discrete) algebraically closed field $k$ over $\F_q$. 
\begin{theorem}[{cf.\ \Cref{sec:vect-bundl-absol-1-classification-of-absolute-vector-bundle}}]
  \label{sec:introduction-main-theorem-introduction}
  The natural functor
  \[
   \mathcal{E}_k(-)\colon \Isoc_k\to \Bun_\FF(\Spd(k))
  \]
  from the category of $E$-isocrystals over $k$ to $\Bun_\FF(\Spd(k))$ is an equivalence of categories. 
\end{theorem}

One might call objects in $\Bun_\FF(\Spd(k))$ ``vector bundles on the absolute Fargues--Fontaine curve'' although we do not literally construct a relative Fargues--Fontaine curve ``$X_{E,\Spd(k)}$'' and define $\Bun_\FF(\Spd(k))$ as vector bundles on ``$X_{E,\Spd(k)}$''.

Let us now recall the definition of $E$-isocrystals over $k$. Let
\[
  W_{\mathcal{O}_E}(-):=\mathcal{O}_E\widehat{\otimes}_{W(\F_q}W(-)
\]
be the functor of $\mathcal{O}_E$-typical Witt vectors on perfect $\F_q$-algebras, and $\varphi\colon W_{\calO_E}(-)\to W_{\calO_E}(-)$ the natural lift of the $q$-Frobenius $x\mapsto x^q$.
Then for any perfect $\F_q$-algebra $R$ the category $\Isoc_R$ is the category of pairs $(D,\varphi_D)$ with $D$ a finite projective $W_{\calO_E}(R)[1/p]$-module and $\varphi_D\colon \varphi^\ast D\to D$ an isomorphism and morphisms respecting the $\varphi_D$.

Given any perfectoid space $S$ over $\Spd(R)$ there exists a natural functor
\[
  \Isoc_R\to \Bun_\FF(S)
\]
constructed as follows:
Let $Y_{E,S}$ be the usual $\Z$-covering of $X_{E,S}=Y_{E,S}/\varphi^\Z$ from \cite[Definition II.1.15]{fargues2021geometrization}. Then to $(D,\varphi_D)\in \Isoc_R$ can be associated the vector bundle on $X_{E,S}$ obtained by descent of the trivial bundle $D\otimes_{W_{\calO_E}(R)[1/p]} \calO_{Y_{E,S}}$ along the $\varphi$-semilinear automorphism $\varphi_D\otimes\varphi$.

By naturality of this construction with respect to pullback along morphisms $S^\prime\to S$ of perfectoid spaces over $\Spd(R)$, we obtain the functor
\[
  \mathcal{E}_R(-)\colon \Isoc_R\to \Bun_\FF(\Spd(R))
\]
mentioned in \Cref{sec:introduction-main-theorem-introduction} if $R=k$.

If $G/E$ is a reductive group, then we will deduce from \Cref{sec:introduction-main-theorem-introduction} a similar description of ``$G$-bundles on the absolute Fargues--Fontaine curve'' as alluded to in the title, cf.\ \Cref{sec:texorpdfstr-bundl-ab-1-classification-of-absolute-bundles}.

Another case that we consider is that of $Y=[\Spd(k)/\underline{H}]$ being the classifying stack of some locally profinite group $H$. In this case $\Bun_\FF(Y)$ is equivalent to continuous representations of $H$ on isocrystals, cf.\ \Cref{sec:vect-bundl-class-1-classification-for-classifying-stack}.

By \cite[Theorem 2.1]{anschutz2021fourier} the functor sending $S\in \Perf_{\F_q}$ to the $\infty$-category of perfect complexes on $X_{E,S}$ is a $v$-stack of $\infty$-categories, which we call $\mathcal{P}erf_\FF$. We make the following conjecture.

\begin{conjecture}
  \label{sec:introduction-conjecture-on-values-on-perfect-schemes}
  For any perfect $\F_q$-algebra $R$ the category $\mathcal{P}erf_\FF(\Spd(R))$ is equivalent to the category of perfect complexes of isocrystals over $R$.\footnote{By \cite[Proposition 2.7]{anschutz2021fourier} perfect complexes of isocrystals over $R$ identify with perfect complexes over $W_{\calO_E}(R)[1/p]$ with a $\varphi$-semilinear isomorphism.}
\end{conjecture}

We note that \Cref{sec:introduction-conjecture-on-values-on-perfect-schemes} includes a comparison of cohomology, which is already non-trivial if $R=k$. Namely, it implies that $R\Hom_{\mathcal{P}erf_\FF(\Spd(k))}(\mathcal{O},\mathcal{O}(-1))=0$.

Let us now go through the contents of the sections of this paper. In \Cref{sec:pro-etale-v} we describe the pro-\'etale site of $\Spd(k)$. Although the outcome is as expected the pro-\'etale site of $\Spec(k)$, the proof uses some difficult results from \cite{scholze_etale_cohomology_of_diamonds}. In \Cref{sec:vect-bundl-absol} we prove \Cref{sec:introduction-main-theorem-introduction}. In \Cref{sec:vect-bundl-class} we apply \Cref{sec:introduction-main-theorem-introduction} to describe the category $\Bun_\FF([\Spd(k)/\underline{H}])$ if $H$ is a locally profinite group. If $G$ is a reductive group over $E$, then in \Cref{sec:texorpdfstr-bundl-ab} we describe $\Bun_G(\Spd(k))$ for the $v$-stack $\Bun_G$ of $G$-bundles on the Fargues--Fontaine curve.

\subsection*{Acknowledgments}
The author wants to express his gratitude to L.\ Fargues for the discussions culminating in the question on absolute vector bundles, which is answered in this paper. Moreover, he thanks Ian Gleason and Kiran S.\ Kedlaya for discussions related to this article. This paper is a completely rewritten version of a preprint from 2016. Namely, back in 2016 the foundations of the theory were not as settled as they are now thanks to \cite{fargues2021geometrization} and \cite{scholze_etale_cohomology_of_diamonds}. This was the authors motivation to revisit the paper.

\section{Pro-\'etale {$v$}-sheaves over {$\Spd(k)$}}
\label{sec:pro-etale-v}

Let $k$ be an algebraically closed field extension of $\F_q$. In this section we analyze small $v$-sheaves $Y$ with a quasi-pro-\'etale map $Y\to \Spd(k)$.

Let $\Spec(k)_{\proet}$ be the pro-\'etale site of $\Spec(k)$ in the sense of \cite[Definition 4.1.1]{bhatt_scholze_the_pro_etale_topology_for_schemes}, or equivalently the site of locally profinite sets with topology induced by surjections of profinite sets. Explicitly, if $T$ is a locally profinite set, then the topological space
\[
  T\cong T\times |\Spec(k)|
\]
is a scheme, when equipped with the pullback of the structure sheaf from $\Spec(k)$.\footnote{\cite[Lemma 2.2.8 ]{bhatt_scholze_the_pro_etale_topology_for_schemes} deals with a similar construction.} We denote this scheme over $\Spec(k)$ by $T\times \Spec(k)$. Clearly, the functor $T\to T\times \Spec(k)$ commutes with limits, and defines the aforementioned equivalence of sites.

Let $S$ be any perfectoid space over $\Spd(k)$, and $T$ a locally profinite set. Then we can define the perfectoid space
\[
  T\times S
\]
by equipping the topological product $T\times |S|$ with the pullback of the structure sheaf on $S$. Clearly, the morphism $T\times S\to S$ is quasi-pro-\'etale in the sense of \cite[Definition 10.1]{scholze_etale_cohomology_of_diamonds}.

Let $S_\proet$ be the pro-\'etale site of $S$ in the sense of \cite[Definition 8.1]{scholze_etale_cohomology_of_diamonds}. Then we get a functor
\[
  Z\in \Spec(k)_\proet\mapsto Z^\diamond\times_{\Spd(k)}S\in S_\proet
\]
with $(-)^\diamond$ the functor from perfect schemes to small $v$-sheaves introduced \cite[Section 18.3]{scholze2020berkeley}. More explicitly, $Z=T\times \Spec(k)$ is sent to $T\times S\cong $ for $T$ a locally profinite set.

By descent we obtain a functor $T\times \Spec(k)\mapsto T\times \Spd(k)$ from $\Spec(k)_\proet$ to the category $\Spd(k)_{\proet}$ of small $v$-sheaves $Y$ with a quasi-pro-\'etale morphism to $\Spd(k)$. Via surjections of $v$-sheaves, the category $\Spd(k)_\proet$ is a site.

\begin{proposition}
  \label{sec:pro-etale-v-1-quasi-pro-etale-v-sheaves-over-spd-k}
  The functor $\Spec(k)_\proet\to \Spd(k)_\proet$ is an equivalence of sites.
\end{proposition}
\begin{proof}
  Let $C/k$ be a non-archimedean, algebraically closed extension of $k$ and $S:=\Spa(C,\mathcal{O}_C)$. By \cite[Proposition 9.6]{scholze_etale_cohomology_of_diamonds} each small $v$-sheaf $Y^\prime$ with a quasi-pro-\'etale map $Y^\prime\to S$ is represented by a perfectoid space (note that the definition of quasi-pro-\'etale in \cite[Definition 10.1]{scholze_etale_cohomology_of_diamonds} includes the assumption of local separatedness).
  By $v$-descent for quasi-pro-\'etale maps, it suffices therefore to show that $\Spec(k)_\proet$ is equivalent to the category of descent data for quasi-pro-\'etale maps for the groupoid $S\times_{\Spd(k)}S\rightrightarrows S$. Now, \cite[Theorem 19.5]{scholze_etale_cohomology_of_diamonds} implies that $S\times_{\Spd(k)} S$ is connected. This implies that for a locally profinite set $T$ the restriction along the diagonal $S\to S\times_{\Spd(k)} S$ induces a bijection
  \[
    \Hom_{\cont}(S\times_{\Spd(k)}S,T)\to \Hom_{\cont}(S,T).
  \]
  This in turn implies that on $T\times S$, there can only exist the trivial descent datum with descent $T\times \Spd(k)$. As $\Spec(k)_\proet$ is equivalent to $S_\proet$ this finishes the proof.
\end{proof}

\begin{corollary}
  \label{sec:pro-etale-v-2-quasi-pro-etale-torsors-on-spd-k}
  Let $H$ be any topological group and $\underline{H}:=\Hom_\cont(-,H)$ its associated constant sheaf on the $v$-site $\Perf_{\F_q}$. Then any quasi-pro-\'etale $H$-torsor on $\Spd(k)$ is trivial. 
\end{corollary}
\begin{proof}
  By \Cref{sec:pro-etale-v-1-quasi-pro-etale-v-sheaves-over-spd-k} we can conclude that each quasi-pro-\'etale morphism $Y\to \Spd(k)$ with $Y$ a small $v$-sheaf, has a section because the same is true for $\Spec(k)_\proet$.
\end{proof}

\section{Vector bundles on the absolute Fargues--Fontaine curve}
\label{sec:vect-bundl-absol}

Let $k$ be an algebraically closed field extension of $\F_q$. In this section we want to prove \Cref{sec:introduction-main-theorem-introduction}, i.e., that the functor
\[
  \mathcal{E}_k(-)\colon \Isoc_k\to \Bun_\FF(\Spd(k))
\]
is an equivalence.

We start with the following lemma.

\begin{lemma}
  \label{sec:vect-bundl-absol-1-absolute-bundles-is-abelian}
  The category $\Bun_\FF(\Spd(k))$ is abelian.
\end{lemma}
\begin{proof}
  Clearly, $\Bun_\FF(\Spd(k))$ is additive because $\Bun_\FF(S)$ is additive for any perfectoid space $S$ over $\Spd(k)$.
  Set $S:=\Spa(C,\calO_C)$ with $C/k$ a non-archimedean, algebraically closed extension with residue field $k$. Let $f\colon S\to \Spd(k)$ be the natural morphism, which induces the exact pullback functor 
\[
f^\ast\colon \Bun_\FF(\Spd(k))\to \Bun_{\FF}(S)\cong \Bun(X_{E,S})\cong \Bun(X_{E,S}^\sch).
\] 
Here we used GAGA for the curve (\cite[Section II.2.3]{fargues2021geometrization}) to identify vector bundles on $X_{E,S}$ with vector bundles on the schematic Fargues--Fontaine curve $X_{E,S}^\sch$. As $S\to \Spd(k)$ is a $v$-cover an absolute vector bundle, i.e., an object of $\Bun_\FF(\Spd(k))$, is equivalently given by a vector bundle $\cal{E}$ on $X_{E,F}$ (equivalently on $X^{\sch}_{E,S}$) together with a descent datum 
\[
\sigma_\cal{E}\colon \pr_1^\ast\cal{E}\cong\pr_2^\ast \cal{E}
\] 
on $X_{E,S\times_{\Spd(k)} S}$ (note that no schematic Fargues--Fontaine curve is defined for $S\times_{\Spd(k)} S$ as $S\times_{\Spd(k)}S$ is not affinoid).
Let $A$ be the group of automorphisms of $S$ over $\Spd(k)$. Pulling back the descent datum along the morphism  
\[
S\times A=\coprod\limits_{a\in A} S\times\{a\}\to S\times_{\Spd(k)}S, \quad (x,g)\mapsto (x,gx)
\]
shows that $\mathcal{E}$ is naturally equipped with a semilinear action of $A$. Recall now that $X_{E,S}^\sch$ is a connected noetherian Dedekind scheme (\cite[Proposition II.2.9]{fargues2021geometrization}), and hence admits a good theory of coherent sheaves. Let $g\colon (\cal{E},\sigma_\cal{E})\to (\cal{F},\sigma_\cal{F})$ be a morphism of descent data, with $\mathcal{E}, \mathcal{F}$ vector bundles on $X^\sch_{E,S}$. We claim that the kernel and cokernel of $g$ are locally free, which then implies that they admit compatible descent data and that $\Bun_\FF(\Spd(k))$ is abelian. To show torsion freeness of the kernel and cokernel, we may enlarge $C$ and in particular assume that it is a field of Mal'cev-Neumann series over $k$, cf.\ \cite[Section 3]{poonen1993maximally}. Then note that the morphism $g\colon \mathcal{E}\to \mathcal{F}$ is $A$-equivariant for the natural semilinear $A$-actions on $\mathcal{E}, \mathcal{F}$. Thus, the torsion in the cokernel and the torsion of the kernel of $g$ aquire natural semilinear $A$-actions.
By \Cref{sec:vect-bundl-absol-2-orbits-of-a-on-schematic-ff-curve} they are therefore trivial as they are supported at an $A$-equivariant, finite set of closed points of $X^\sch_{E,S}$. This finishes the proof.
\end{proof}

Let us recall the definition of the field of Mal'cev-Neumann series $k((t^\R))$ over $k$, cf.\ \cite[Section 3]{poonen1993maximally}. It is the ring of formal power series
  \[
    a=\sum\limits_{x\in \R} a_xt^x
  \]
  with coefficients $a_x\in k$, such that the support
  \[
    \mathrm{supp}(a):=\{x\in \R \ |\ a_x\neq 0\}
  \]
  is a well-ordered subset of $\R$.

The multiplication
\[
  (\sum\limits_{x\in \R}a_xt^x)(\sum\limits_{x\in \R} b_xt^x):=\sum\limits_{x\in \R}(\sum\limits_{y+z=x}a_yb_z)t^x
\]
is well-defined due to the support condition and makes $k((t^\R))$ into a field, cf.\ \cite[Section 3]{poonen1993maximally}.
The natural valuation
\[
  \nu\colon k((t^\R))\to \R\cup \{\infty\},\ a\mapsto \mathrm{inf}\{\mathrm{supp}(a)\}
\]
shows that $k((t^\R))$ is naturally a non-archimedean field, which is even algebraically closed, cf.\ \cite[Theorem 1]{poonen1993maximally}. 

The construction of $k((t^\R))$ works with $\R$ replaced by any ordered group, and is in fact functorial in automorphisms of such. In particular, the multiplication by some $x\in \R_{>0}$ induces a continuous automorphism
\[
  \sigma_x\colon k((t^\R))\to k((t^\R)),
\]
which sends $t$ to $t^x$.

\begin{lemma}
  \label{sec:vect-bundl-absol-2-orbits-of-a-on-schematic-ff-curve}
  Set $C:=k((t^\R))$, $S:=\Spa(C,\calO_C)$ and $A:=\mathrm{Aut}_{\Spd(k)}(S)$. Then the action of the group $A$ on the closed points of $X^{\sch}_{E,S}$ has no non-empty finite orbit.
\end{lemma}
\begin{proof}
  Assume that $\{x_1,\ldots, x_n\}\subseteq X^{\sch}_{E,S}$ is a non-empty finite $A$-orbit.
  As the morphism $X_{E,S}\to X^{\sch}_{E,S}$ induces a bijection on classical points (\cite[Proposition II.2.9]{fargues2021geometrization}), we can view $\{x_1,\ldots, x_n\}$ as a subset of $X_{E,S}$. Let
  \[
    f\colon Y_{E,S}\to X_{E,S}=Y_{E,S}/\varphi^{\Z}
  \]
  be the natural covering, and let $y_i\in Y_{E,S}$ be a classical point mapping to $x_i$, cf.\ \cite[Definition/Proposition  II.1.22]{fargues2021geometrization}. Then by the same reference
  \[
    f^{-1}(\{x_1,\ldots, x_n\})=\{\varphi^i(y_1),\ldots,\varphi^i(y_n)\ |\ i\in \Z \}.
    \]
  By \cite[Proposition II.1.8]{fargues2021geometrization} there exists non-units $a_1,\ldots, a_n\in \mathcal{O}_C$ such that $y_i$ is the vanishing locus of $\pi-[a_i]$, with $\pi\in \calO_E$ a fixed uniformizer.
  Let $\nu\colon C\to \R$ be the $t$-adic valuation introduced above.
  The set
  \[
    M:=\{\nu(\varphi^i(a_j))\ |\ i\in \Z, j=1,\ldots, n\}\subseteq \R
  \]
  is then countable. This implies that the subgroup $\langle M\rangle\subseteq \R$ generated by $M$ is countable. Hence, there exists some $x\in \R_{> 0}\setminus \langle M\rangle$. Then $x+m\notin M$ for all $m\in M$. Consider now the automorphism $\sigma_x\in A$ sending $t$ to $t^x$ from above. Clearly,
  \[
    \nu(\sigma_x(a))=x+\nu(a)
  \]
  for any $a\in C$. This implies that $\sigma_x(y_1)\notin f^{-1}(\{x_1,\ldots, x_n\})$. Indeed, assume that $\sigma_x(y_1)=\varphi^{m}(y_j)$ for some $m\in \Z$ and $j=1,\ldots, n$, and let $\theta\colon W_{\mathcal{O}_E}(\mathcal{O}_C)\to k(\varphi^m(y_j))$ be the natural map with kernel $(\pi-[a_j^{p^m}])$. Then $\sigma_x(y_1)=\varphi^m(y_j)$ implies
  \[
    \theta(\pi-[\sigma_x(a_1)])=0
  \]
  and thus $\theta([\sigma_x(a_1)])=\sigma_x(a_1)^\sharp=(a_{j}^{p^m})^\sharp=\theta([a_{j}^{p^m}])$ for the $\sharp$-map from the natural identification $C\cong k(\varphi^m(y_j))^\flat$. But this is impossible because $\sigma_x(a_1)$ and $a_j^{p^m}$ have different valuations. This finishes the proof.
\end{proof}

Let $C/k$ be a non-archimedean, algebraically closed extension and $S:=\Spa(C,\mathcal{O}_C)$.
By \cite[Section II.2.4]{fargues2021geometrization} each vector bundle on $X_{E,S}$ has its rank and degree, and these yield a Harder-Narasimhan formalism in the sense of \cite[Chapitre 5.5]{fargues_fontaine_courbes_et_fibres_vectoriels_en_theorie_de_hodge_asterisque}.
By descent, we can then define the rank and degree for the absolute vector bundles, i.e., the objects in $\Bun_\FF(\Spd(k))$ as well. Using the identity functor as a ``generic fiber functor'' (as permitted by \Cref{sec:vect-bundl-absol-1-absolute-bundles-is-abelian}) the degree and rank function satisfy the conditions of an abstract Harder-Narasimhan formalism on $\Bun_{\FF}(\Spd(k))$. In particular, each absolute vector bundle has its Harder-Narasimhan filtration with associated gradeds given by semi-stable objects.

We note that if $\mathcal{E}\in \Bun_\FF(\Spd(k))$, then the pullback of its HN-filtration to $S$ is the HN-filtration of the pullback $\mathcal{E}_S\in \Bun_\FF(S)$. Indeed, over $S\times_{\Spd(k)} S$ there exists at most one HN-filtration of an object in $\Bun_\FF(S\times_{\Spd(k)}S)$, and hence the pullbacks of the HN-filtration of $\mathcal{E}_S$ along the projections $S\times_{\Spd(k)}S \to S$ must agree. In particular, the HN-filtration descents to a filtration of $\mathcal{E}$. But the graded pieces for this filtrations are semistable (as they are after pullback to $S$) and hence this filtration must be the HN-filtration of $\mathcal{E}$.

We now analyze semistable objects in $\Bun_{\FF}(\Spd(k))$.

\begin{lemma}
  \label{sec:vect-bundl-absol-1-semistable-bundles-and-morphisms-between-them}
  Let $\lambda,\mu\in \Q$, and $\mathcal{E},\mathcal{F}\in \Bun_\FF(\Spd(k))$ be semistable of slopes $\lambda$ resp.\ $\mu$.
  \begin{enumerate}
  \item $\mathcal{E}$ is isomorphic to a direct sum of $\mathcal{O}(\lambda)$.
    \item If $\lambda\neq \mu$, then $\Hom(\mathcal{E},\mathcal{F})=0$.
  \end{enumerate}
\end{lemma}
Here, $\mathcal{O}(\lambda)=\mathcal{E}_k(D_\lambda)$ is the vector bundle associated with the isoclinic isocrystal $D_\lambda$ of slope $-\lambda$, cf.\ \cite[Section II.2.]{fargues2021geometrization}.
\begin{proof}
  By \cite[Theorem III.4.5]{fargues2021geometrization} the torsor of isomorphisms of $\mathcal{E}$ to the direct sum of $\mathcal{O}(\lambda)$ (of rank $\mathrm{rk}(\mathcal{E})$) is representable by some pro-\'etale $v$-sheaf over $\Spd(k)$. By \Cref{sec:pro-etale-v-1-quasi-pro-etale-v-sheaves-over-spd-k} this torsor splits, showing the first claim.
  If $\lambda>\mu$, then the second statement is a general consequence of the HN-formalism, cf.\ \cite[Section 5.5.]{fargues_fontaine_courbes_et_fibres_vectoriels_en_theorie_de_hodge_asterisque}.
  Thus, let us assume that $\lambda<\mu$. By (1) we know that $\mathcal{E}, \mathcal{F}$ are direct sums of $\mathcal{O}(\lambda)$'s resp.\ $\mathcal{O}(\mu)$'s. Considering $\mathcal{E}^\vee\otimes \mathcal{F}$ then reduces to the case that $\mathcal{E}\cong \mathcal{O}$ and $\mu>0$. Considering pushforwards along $g_S\colon X_{E^\prime,S}\to X_{E,S}$ for a finite, separable extension $E^\prime|E$, then reduces to the case that $\mu\in \Z$. More precisely, for any perfectoid space $S^\prime$ over $\Spd(k)$ the pullback $\mathcal{F}_{S^\prime}$ embeds via the unit of the adjunction into $g_{S^\prime, \ast}g^{\ast}_{S^\prime} \mathcal{F}_{S^\prime}$ and this embedding commutes with base change in $S^\prime$. Hence, it descents to $\Spd(k)$. Now,
  \[
    \Hom_{X_{E,S^\prime}}(\mathcal{O},g_{S^\prime,\ast}g^\ast_{S^\prime} \mathcal{F}_{S^\prime})\cong \Hom_{X_{E^\prime,S^\prime}}(\mathcal{O},g^\ast_{S^\prime} \mathcal{F}_{S^\prime}),
  \]
  again compatible with base change, and this allows to replace $E$ by $E^\prime$ and $\mathcal{F}$ by $g^\ast_{S^\prime} \mathcal{F}_{S^\prime}$.
  If $\mu\in \Z_{>0}$, then each $\Hom_{\Bun_\FF(\Spd(k))}(\mathcal{O},\mathcal{O}(\mu))=0$ as each non-zero morphism must have a torsion cokernel after pullback to some algebraically closed, non-archimedean field of Mal'cev-Neumann series and thus is surjective by \ref{sec:vect-bundl-absol-2-orbits-of-a-on-schematic-ff-curve}. This finishes the proof. 
\end{proof}

We can now prove the desired classification of absolute vector bundles.

\begin{theorem}
  \label{sec:vect-bundl-absol-1-classification-of-absolute-vector-bundle}
  Let $k$ be an algebraically closed field over $\F_q$. Then the functor
  \[
    \mathcal{E}_k(-)\colon \Isoc_k\to \Bun_\FF(\Spd(k))
  \]
  is an equivalence.
\end{theorem}
\begin{proof}
  Fully faithfulness follows from \Cref{sec:vect-bundl-absol-1-semistable-bundles-and-morphisms-between-them} and \cite[Proposition II.2.5]{fargues2021geometrization}. Indeed, by the Dieudonn\'e-Manin classification the category $\Isoc_k$ is semisimple, and for each simple object $D_\lambda$ with slope $\lambda\in \Q$, \cite[Proposition II.2.5]{fargues2021geometrization} implies that $\mathrm{End}(D_\lambda)\cong \mathrm{End}(\mathcal{E}_k(D_\lambda))$ (in fact, this statement holds even for $\Spd(k)$ replaced by $\Spa(C, \calO_C)$ with $C/k$ non-archimedean, algebraically closed).
  Now \Cref{sec:vect-bundl-absol-1-semistable-bundles-and-morphisms-between-them} implies that there are no non-zero morphisms between semistable objects for different slopes, which shows then fully faithfulness.
  
  To show essential surjectivity we argue by induction on the rank. The case of rank $0$ is trivial. Let $\mathcal{E}\in \Bun_\FF(\Spd(k))$.
  We first assume that $\mathcal{E}$ is semistable of slope $\lambda$. In this case, we claim that each non-zero subobject $\mathcal{F}\subsetneq \mathcal{E}$ has slope $\lambda$. By the proof of \ref{sec:vect-bundl-absol-1-absolute-bundles-is-abelian} we must have $\mathrm{rk}(\mathcal{F})<\mathrm{rk}(\mathcal{E})$. By induction $\mathcal{F}$ can thus be assumed to lie in the essential image of $\mathcal{E}_k(-)\colon \Isoc_k\to \Bun_\FF(\Spd(k))$. In particular, $\mathcal{F}$ is a sum of $\mathcal{O}(\mu)$'s for some (different) $\mu\in \Q$. As $\mathcal{F}\to \mathcal{E}$ is injective, we can deduce from \Cref{sec:vect-bundl-absol-1-semistable-bundles-and-morphisms-between-them} that in fact $\mathcal{F}$ must be of slope $\lambda$.

  Now, we claim that for each $\mathcal{E}\in \Bun_\FF(\Spd(k))$ its HN-filtration is canonically split. For this we want to use the trick of considering the negative degree function from \cite[5.5.2.3]{fargues_fontaine_courbes_et_fibres_vectoriels_en_theorie_de_hodge_asterisque} to get an opposite HN-filtration, which yields a canonical splitting. However, a priori the semistable objects in $\Bun_\FF(\Spd(k))$ for the usual HN-formalism need not be semistable for the opposite HN-formalism for $(-\mathrm{deg},\mathrm{rk})$. However, we checked above that each semistable object of slope $\lambda$ only has only non-trivial subobjects of slope $\lambda$. In particular, they remain semistable for the opposite HN-filtration. This implies that the usual HN-filtration splits the opposite one, and hence the claim.

  Using induction on the rank and \Cref{sec:vect-bundl-absol-1-semistable-bundles-and-morphisms-between-them} the essential surjectivity follows.
\end{proof}

\section{Vector bundles for classifying stacks}
\label{sec:vect-bundl-class}

Let $k$ be an algebraically closed extension of $\F_q$, and let $H$ be a locally profinite group. We let $\underline{H}:=\Hom_\cont(-,H)$ be the associated small $v$-sheaf on $\Perf_{\F_q}$. Let $[\Spd(k)/\underline{H}]$ be the classifying stack of $\underline{H}$ over $\Spd(k)$. In this section, we want to identify the category
\[
  \Bun_\FF([\Spd(k)/\underline{H}]).
\]
By $v$-descent, this category identifies with descent data for the groupoid
\[
  \underline{H}\times \Spd(k)\rightrightarrows \Spd(k).
\]
Note that $\underline{H}\times \Spd(k)\cong (H\times \Spec(k))^\diamond$ in the notation of \Cref{sec:pro-etale-v}. In particular, $\underline{H}\times \Spd(k)$ is the $v$-sheaf associated to some perfect scheme over $\Spec(k)$.

We need the following lemma.
\begin{lemma}
  \label{sec:vect-bundl-class-1-fully-faithful-for-pullback}
  Let $f\colon Z\to \Spec(k)$ be a perfect scheme over $\Spec(k)$, let $f^\diamond\colon Z^\diamond\to \Spd(k)$ be the natural map and let $D,D^\prime\in \Isoc_k$. Then the map
  \[
    \Hom_{\Isoc_Z}(f^\ast D, f^\ast D^\prime)\to \Hom_{\Bun_\FF(\Spd(k))}(f^{\diamond, \ast}\mathcal{E}_k(D), f^{\diamond, \ast}\mathcal{E}_k(D^\prime))
  \]
  is an isomorphism.
\end{lemma}
Here, the left hand side denotes morphisms of $E$-isocrystals over $Z$.
\begin{proof}
  By the Dieudonn\'e-Manin classification, we may assume that $D, D^\prime$ are simple of slope $\lambda$ resp.\ $\lambda^\prime$. We may even reduce to the case that $\lambda=0$, i.e., $D=W_{\calO_E}(k)[1/p]$ with its standard Frobenius. As the statement is local on $Z$ we may assume that $Z=\Spec(R)$ is affine.
  If $\lambda^\prime\neq \lambda$ we claim that both sides vanish. As $R$ is reduced, it embeds into a product $\prod\limits_{i\in I}k_i$ of algebraically closed fields. Then $W_{\calO_E}(R)[1/p]$ embeds into $\prod\limits_{i\in I}(W_{\calO_E}(k_i)[1/p])$ and from the case of an algebraically closed field we can deduce by the Dieudonn\'e-Manin classification that the left hand side vanishes.

  Let us show the vanishing of the right hand side. If $\lambda^\prime<0$, then the vanishing follows from \cite[Proposition II.2.5.(1)]{fargues2021geometrization} by descent. Thus, assume that $\lambda^\prime>0$. Rephrased using absolute Banach-Colmez spaces, we have to show that each morphism
  \[
    g\colon \Spd(R)\to \mathcal{BC}(\mathcal{O}(\lambda^\prime))
  \]
  over $\Spd(k)$ factors over the zero section if $\lambda^\prime>0$. This may be checked after base change along some $S:=\Spa(C, \mathcal{O}_C)\to \Spd(k)$ with $C$ non-archimedean and algebraically closed. Over $S$, it suffices to test this on geometric points $S^\prime=\Spa(C^\prime, C^{\prime, +})$ over $S$. Now each morphism $R\to C^\prime$ factors over some algebraic closure of some residue field of $R$. Hence, we may replace $R$ by such a field, and in this case the desired vanishing follows from \Cref{sec:vect-bundl-absol-1-semistable-bundles-and-morphisms-between-them}.
  To finish the proof, we have to understand the case that $\lambda^\prime=0$. Then the left hand side identifies with $\underline{E}(\Spec(R))$ by Artin-Schreier theory, while the right hand side with $\underline{E}(\Spd(R))$. By \cite[Proposition 18.2.1]{scholze2020berkeley} the open and closed decompositions of $\Spec(R)$ and $\Spd(R)$ agree (this follows by considering homomorphisms to $\Spec(\F_p\times \F_p)$ whose associated $v$-sheaf is $\Spd(\F_p)\coprod \Spd(\F_p)$). In particular, $\pi_0(\Spec(R))$ and $\pi_0(\Spd(R))$ are homeomorphic. This implies the claim and the proof is finished.  
\end{proof}

We can now identify the category $\Bun_\FF([\Spd(k)/\underline{H}])$. By \cite[Proposition 5.10]{ivanov2020ind} the fibered category $Z\mapsto \Isoc_Z$ is a $v$-stack on perfect $\F_q$-schemes (even an arc-stack) and hence we can talk about $\Isoc_{[\Spec(k)/\underline{H}]}$, which more concretely identifies with the category of descent data for the groupoid
\[
  \underline{H}\times \Spec(k)\rightrightarrows \Spec(k).
\]

\begin{theorem}
  \label{sec:vect-bundl-class-1-classification-for-classifying-stack}
  The natural functor $\Isoc_{[\Spec(k)/\underline{H}]}\to \Bun_\FF(\Spd(k)/\underline{H}])$ is an equivalence.
\end{theorem}
\begin{proof}
  This follows from \Cref{sec:vect-bundl-class-1-fully-faithful-for-pullback} and \Cref{sec:vect-bundl-absol-1-classification-of-absolute-vector-bundle}. Namely, both imply that each descent datum on some object $\mathcal{E}\in \Bun_{\FF}(\Spd(k))\cong \Isoc_k$ is already defined on the corresponding isocrystal.
\end{proof}

Assume now that $H$ is profinite and let $C(H,k)$ be the ring of continuous functions on $H$ with values in $k$.
Then
\[
  \underline{H}\times \Spec(k)\cong \Spec(C(H,k)).
\]
Moreover, the natural maps
\[
C(H, \Z_p)\widehat{\otimes}_{\Z_p}W_{\calO_E}(k)\to C(H,W_{\calO_E}(k))\to W_{\calO_E}(C(H,k))
\]
are isomorphisms, where the topology on $W_{\calO_E}(k)$ is the $p$-adic topology and the tensor product is $p$-adically complete. In fact, if $\pi\in \mathcal{O}_E$ is a uniformizer the three objects are $\pi$-torsion free, $\pi$-complete lifts of the perfect $k$-algebra $C(H,k)$.
Inverting $p$ we get isomorphisms
\[
  C(H,\Q_p)\widehat{\otimes}_{\Q_p}W_{\calO_E}(k)[1/p]\cong C(H,W_{\calO_E}(k)[1/p])\cong W_{\calO_E}(C(H,k))[1/p].
\]
Objects in $\Isoc_{[\Spec(k)/\underline{H}]}$ identify therefore with isocrytals $D\in \Isoc_{k}$ together with a coaction of the Hopf algebra $C(H,\Q_p)$. Unraveling the coaction, this yields precisely an isocrystal $D$ with a continuous action of $H$ on $D$, when $D$ is equipped with its natural topology as a finite dimensional $W_{\calO_E}(k)[1/p]$-vector space. Hence, objects in $\Isoc_{[\Spd(k)/\underline{H}]}$ deserve to be called ``isocrystals with continuous $H$-action'', even if $H$ is only locally profinite.

\section{\texorpdfstring{$G$-bundles on the absolute Fargues--Fontaine curve}{G-bundles on the absolute Fargues--Fontaine curve}}
\label{sec:texorpdfstr-bundl-ab}

Let $G$ be a reductive group over $E$. In this section we want to generalize \Cref{sec:vect-bundl-absol-1-classification-of-absolute-vector-bundle} from $\GL_n$-bundles to arbitrary $G$-bundles on the absolute Fargues--Fontaine curve.

First we recall the definition of a $G$-bundle that we use, cf.\ \cite[Theorem 19.5.2]{scholze2020berkeley}, \cite[Definition/Proposition III.1.1]{fargues2021geometrization}.

\begin{definition}
\label{definition-g-bundle}
Let $S/\F_q$ be a perfectoid space. A $G$-bundle (or $G$-torsor) over $X_{E,S}$ is an exact, $E$-linear tensor functor
\[
\mathrm{Rep}_E(G)\to \Bun(X_{E,S})
\]
from the category of representations of $G$ on finite-dimensional $E$-vector spaces to the category of vector bundles on $X_{E,S}$.
\end{definition}

For a perfectoid space $S/\F_q$ we denote by $\Bun_G(S)$ the groupoid of $G$-bundles on $X_{E,S}$ and by $\Bun_G$ the fibered category associating to each $S$ the groupoid $\Bun_G(S)$. As in the case of vector bundles the fibered category $\Bun_G$ is a small $v$-stack, cf.\ \cite[Section III.1]{fargues2021geometrization}

Let $k$ be an algebraically closed field over $\F_q$. As for the case of $\Bun_\FF$, we are interested in the category
\[
  \Bun_G(\Spd(k))
\]
of ``$G$-bundles on the absolute Fargues--Fontaine curve''.

To ease notation, we set
\[
  L:=W_{\calO_E}(k)[1/p],
\]
which carries its natural Frobenius $\varphi\colon L\to L$.

The following gadget will be shown to be equivalent to $\Bun_G(\Spd(k))$.

\begin{definition}
\label{definition-kottwitz-category}
The Kottwitz category $\mathcal{B}(G)$ for $G$ is defined as follows. 
It has as objects the elements in
\[
G(L).
\]
For $b,b'\in G(L)$ the set of homomorphisms from $b$ to $b'$ is defined to be the set
\[
\{c\in G(L)\ | \ cb\varphi(c)^{-1}=b'\}.
\]
Finally, composition is defined by multiplication in $G(L)$.
\end{definition}

In other words, the Kottwitz category is the quotient groupoid 
\[
[G(L)/{\varphi-\mathrm{conj.}}]
\]
of $G(L)$ modulo $\varphi$-conjugacy. The set of isomorphism classes of the Kottwitz category, i.e., the quotient set of $G(L)$ modulo $\varphi$-conjugation, is Kottwitz' famous set $B(G)$, cf.\ \cite{Kottwitz1985}, \cite[Definition III.2.1]{fargues2021geometrization}.

By Steinberg's theorem (\cite[Chapitre III.2.3, Th\'er\`eme 1']{Serre2002}) the Kottwitz category is equivalent to the groupoid of exact, $E$-linear tensor functors
\[
  \Rep_EG\to \Isoc_k.
\]
In particular, we can define for any perfectoid space $S$ over $\Spd(k)$ a functor
\[
  \mathcal{T}_S\colon \mathcal{B}(G)\to \Bun_G(S)
\]
by composing an exact, $E$-linear tensor functor $\Rep_E G\to \Isoc_k$ with the exact, $E$-linear tensor functor $\Isoc_k\to \Bun_\FF(S)$ considered in \Cref{sec:introduction-}. As the different functors $\mathcal{T}_S$ are natural with respect to pullback along some morphism $S^\prime\to S$ over $\Spd(k)$, we get a functor
\[
  \mathcal{T}\colon \mathcal{B}(G)\to \Bun_G(\Spd(k)).
\]

As a consequence of \Cref{sec:vect-bundl-absol-1-classification-of-absolute-vector-bundle} we obtain the following.
\begin{theorem}
  \label{sec:texorpdfstr-bundl-ab-1-classification-of-absolute-bundles}
  The functor $\mathcal{T}\colon \mathcal{B}(G)\to \Bun_G(\Spd(k))$ is an equivalence.
\end{theorem}
\begin{proof}
  We know that
  \[
    \mathcal{B}(G)\cong \Fun^{\ex,\otimes}_{E}(\Rep_EG, \Isoc_k)
  \]
  is the category of exact, $E$-linear tensor functors. By \Cref{sec:vect-bundl-absol-1-classification-of-absolute-vector-bundle} the category $\Isoc_k$ identifies with the category of descent data for the groupoid $S\to \Spd(k)$ for any non-empty perfectoid space $S$ over $\Spd(k)$.
  This implies that the category $\Fun^{\ex,\otimes}_{E}(\Rep_E G, \Isoc_k)$ identifies with the category of descent data along $S\to \Spd(k)$ for the $v$-stack $\Fun^{\ex,\otimes}_E(\Rep_E G,\Bun_\FF)$ of exact, $E$-linear tensor functors from $\Rep_EG$ to $\Bun_\FF$. Namely, exactness can be checked $v$-locally and also an $E$-linear resp.\ a tensor structure can be constructed $v$-locally. By \Cref{definition-g-bundle} the $v$-stack $\Fun^{\ex,\otimes}_E(\Rep_E G,\Bun_\FF)$ is equivalent to the $v$-stack $\Bun_G$. This implies that $\mathcal{B}(G)$ is equivalent to $\Bun_G(\Spd(k))$, and in fact the constructed equivalence is given by $\mathcal{T}$. 
\end{proof}

\printbibliography

\end{document}